\documentclass{siamart1116}%
\usepackage{amsmath}
\usepackage{amsfonts}
\usepackage{amssymb}
\usepackage{graphicx}
\usepackage{lipsum}
\usepackage{amsfonts}
\usepackage{graphicx}
\usepackage{epstopdf}
\usepackage{amsopn}
\usepackage{enumerate}
\usepackage{float}
\usepackage{indentfirst}
\usepackage{verbatim}
\usepackage{algorithm}
\usepackage{algpseudocode}
\usepackage{float}
\usepackage{multicol}
\usepackage{listings}
\usepackage{subfigure}
\usepackage{epsfig}
\usepackage{amsfonts}
\usepackage{amssymb}
\usepackage[justification=centering]{caption}
\usepackage{epstopdf}%
\ifpdf
\DeclareGraphicsExtensions{.eps,.pdf,.png,.jpg}
\else
\DeclareGraphicsExtensions{.eps}
\fi
\numberwithin{theorem}{section}
\newcommand{\TheTitle}{A Path Integral Monte Carlo Method based on Feynman-Kac Formula for Electrical Impedance Tomography}
\newcommand{\TheAuthors}{ Yijing Zhou Wei Cai}
\headers{A Path Integral Monte Carlo Method for EIT}{\TheAuthors}

\ifpdf
\hypersetup{
pdftitle={\TheTitle},
pdfauthor={\TheAuthors}
}
\fi
\begin{document}

\title{{A Path Integral Monte Carlo Method based on Feynman-Kac Formula for
Electrical Impedance Tomography}\thanks{Submitted to the editors DATE.
\funding{ W. Cai was supported by  US National Science Foundation (Grant No. DMS-1764187).}}%
}
\author{Yijing Zhou \thanks{Department of Mathematics and Statistics, University of
North Carolina at Charlotte, Charlotte, NC 28223, USA.}
\and Wei Cai\thanks{Department of Mathematics Southern Methodist University, Dallas, TX75275. Corresponding author,
cai@smu.edu.} }
\maketitle

\begin{abstract}
A path integral Monte Carlo method (PIMC) based on Feynman-Kac formula for
mixed boundary conditions of elliptic equations is proposed to solve the
forward problem of electrical impedance tomography (EIT) on the boundary to
obtain electrical potentials. The forward problem is an important part for
iterative algorithms of the inverse problem of EIT, which has attracted
continual interest due to its applications in medical imaging and material
testing of materials. By simulating reflecting Brownian motion with
walk-on-sphere techniques and calculating its corresponding local time, we are
able to obtain accurate voltage-to-current map for the conductivity equation
with mixed boundary conditions for a 3-D spherical object with eight
electrodes. Due to the local property of the PIMC method, the solution of the
map can be done locally for each electrode in a parallel manner.

\end{abstract}




\begin{keywords}
EIT, Reflecting Brownian Motion, boundary local time, Skorohod problem,
Brownian motion, conductivity equation, Feynman-Kac formula, WOS, mixed
boundary problem, boundary element method
\end{keywords}

\begin{AMS}
60G60 
62P30 
78M50 
\end{AMS}

\section{Introduction}

Electrical Impedance Tomography (EIT) is a non-invasive medical imaging
technique in which an image of the electrical properties (conductivity or
permittivity) of part of the body is inferred from surface electrode
measurements. It has the advantages over other techniques like X-rays and
requires no exposure to radioactive materials. Applications include detection
of breast cancer, pulmonary emboli, blood clots, impaired gastric emptying and
etc.. Essentially, through only surface measurements, the internal electric
conductivity and permittivity are identified as an image inside the human
body. For instance, the electric conductivity of malignant tumor, a
high-water-content tissue, is one order higher than that of the normal (fat)
tissue, which allows one to find potential diseases and locations through the
constructed image throughout the body \cite{[Cheney1999]}.

EIT is also a useful tool in other fields such as geophysics, environmental
sciences and nondestructive testing of materials. It is able to locate
underground mineral deposits, detect leaks in underground storage tanks and
monitor flows of injected fluids into the earth for extraction or
environmental cleaning. Moreover, EIT can detect the the corrosion or defects
of construction material and machine parts \cite{[Borcea2002]}
\cite{[Cheney1999]} when invasion testing is not possible or destructive.

For the applications mentioned above, researchers are often faced with the
problem of how to work out the conductivity inside an object when only part of
the boundary measurements are available \cite{[Alessandrini1998]}. It is well
known that this inverse problem is nonlinear, unstable and intrinsically
ill-posed \cite{[Borcea2002]}.

In theory, complete boundary measurements indeed determine the conductivity in
the interior uniquely \cite{[Kohn1984]} \cite{[Sylvester1986]}, however, in
practice only limited number of electrodes and current patterns are available
from measurements. Various numerical algorithms to reconstruct the
conductivity have been proposed and fall into two categories, non-iterative
and interative methods. Noniterative methods were developed based on the
assumption that the conductivity does not differ too much from a constant.
Calder\'{o}n \cite{[Calderon2006]} proved that a map between the conductivity
constant $\gamma$ and a quadratic energy functional $Q_{\gamma}$ in
(\ref{eq6-4}) is injective when $\gamma$ is close to a constant in a
sufficiently small neighborhood, and an approximation formula was also given
to reconstruct the conductivity. The back-projection method of Barber-Brown
\cite{[Santosa1988]} gave a crude approximation to the conductivity increment
$\delta\gamma$ based on inverse of the generalized Radon transform, which
works best for smooth $\delta\gamma$ or $\delta\gamma$ whose singularity is
far from the boundary. For a L-electrode system, Noser algorithm minimizes the
sum of squares error of the voltages on the electrodes by using one step of a
regularized Newton's method. Meanwhile, the iterative methods are devoted to
minimize different regularized least squares functionals such as a
Tikhonov-type regularization \cite{[Somersalo2000]}\cite{[Vauhkonen1998]} and
a total variation formulation \cite{[Borsic2007]} where iterative
gradient-based optimization algorithms are commonly used.

Solving inverse problems with iterative algorithms usually requires a solution
of a forward problem at each iteration numerically, the computation time
accumulates fast for commonly used grid based global methods such as FEM/BEM
methods. In most of the EIT problem, the measurements are only available on
limited number of electrodes, in finding the conductivity profile to match the
measured voltages there, a global solution of the potentials over the whole
object is in fact not needed, therefore a global solution procedure during hte
forward problem incur unneccessary computational cost beyond the electrodes.
With this in mind, in this paper, an alternative local stochastic approach
based on path integral Monte Carlo (PIMC) simulations using Feynman-Kac
formula will be proposed. Due to the nature of the Feynma-Kac formula which
allows the potential solution at any single location including those on the
electrodes, we could dramatically reduce the amount of solutions needed for
each forward problem solution.

The remainder of the paper is organized as follows. In Section 2, the forward
and inverse problems of the EIT problem are introduced. Section 3 includes the
path integral Monte Carlo method for the forward mixed boundary value problem
of the voltage-to-current map for a 8-electrode EIT problem for a spherical
object. Also, in order to validate the accuracy of the stochastic method, we
include a boundary element method to generate reference solutions. Comparison
between the stochastic and deterministic methods shows the accuracy of the
proposed method. Finally conclusions are drawn in Section 4.


\section{Forward and inverse problems in EIT}


\subsection{The forward problem}

In this section, we will first review the forward problem arising from EIT.
The mathematic models for EIT have been developed and compared with the
experimental measurement of voltages on electrodes for a given conductivity
distribution, which is adjusted to fit the measurements. The existing models
are continuum model, gap model, shunt model and complete electrode model.
Among all, the complete electrode model was shown to be capable of predicting
the experimentally measured voltages to within 0.1 percent
\cite{[Somersalo1992]} and the existence and uniqueness of the model has also
been proved .\medskip

\begin{itemize}
\item Complete Electrode Model\medskip
\end{itemize}

Let the domain of the object is denoted as $\Omega,$ embedded within we assume
there is an anomaly $\Omega_{0}\subset\Omega$. The domain $\Omega$ is assumed
to have a smooth boundary with a limited number of electrodes $E_{i}%
,i=1,...,L$ attached to $\partial\Omega$. The conductivity inside $\Omega$ is
given by $\gamma$ and the electric potential for the model will satisfy the
following boundary value problem
\begin{subequations}
\label{eq6-9}%
\begin{align}
\nabla\cdot\gamma\nabla u  &  =0,\quad\text{in \ }\Omega,\label{eq6-9-1}\\
\int_{E_{l}}\gamma\frac{\partial u}{\partial n}dS  &  =J_{l},\quad
l=1,2,...,L,\label{eq6-9-2}\\
\gamma\frac{\partial u}{\partial n}  &  =0,\quad\mbox{off}\ \ \cup_{l=1}%
^{L}E_{l},\label{eq6-9-3}\\
u+z_{l}\gamma\frac{\partial u}{\partial n}  &  =U_{l}\quad\mbox{on}\ \ E_{l}%
,\quad l=1,2,...,L,\label{eq6-9-4}\\
\sum_{l=1}^{L}J_{l}  &  =0,\\
u  &  =u_{0}\text{ \ \ \ \ \ on }\partial\Omega_{0}. \label{tumorBC}%
\end{align}
In (\ref{tumorBC}), we have prescribed a constant potential $u_{0}$ on the
surface of the tumor, which corresponds to modeling the highly conductive
anomaly as a perfect conductor.

Equation (\ref{eq6-9-1}) describes the distribution of electric potential $u$
in the interior of the object where $\gamma$ is the conductivity or inverse of
the resistivity. It was derived from Maxwell equations by neglecting the
time-dependence of the alternating current and assuming the current source
inside the object to be zero \cite{[Somersalo1992]}. Both equations
(\ref{eq6-9-2}) and (\ref{eq6-9-3}) indicates the knowledge of current density
on and off the electrodes on the boundary, respectively. Equation
(\ref{eq6-9-4}) takes account of the electrochemical effect by introducing
$z_{l}$ as the contact impedance or surface impedance which quantitatively
characterizes a thin, highly resistive layer at the contact between the
electrode and the skin, which causes potential jumps according to the Ohms
law. It should be noted that, the regularity of potential $u$ decreases as the
contact impedance approaches zero \cite{[Maire2015]}, which becomes a huge
hindrance to accurate numerical resolution as in practice usually good
contacts with small contact impedance are used.

The usual choice of numerical method to solve the forward problem is grid
based methods such as finite element method, boundary element method. The
global solution will require the solution of a large linear system as fine
meshes ususally are needed around the electrodes. Our approach in this paper
is to develop a local stochastic method, i.e. a path integral Monte Carlo
method based on on Feynman-Kac formula for the mixed boundary value problem in
(\ref{eq6-9-1})-(\ref{eq6-9-4}). This will allows us to only find the unknown
solution on a single electrode only without solving any global matrix system,
thus dramatically reducing the cost of a forward problem solver in an
iterative inverse problem solution of EIT.


\subsection{The inverse problem}

The EIT inverse problem proposed by Calder\'{o}n in 1980 posed the question if
it was possible to determine the heat conductivity of static temperature and
heat flux measurements on the boundary? The inverse problem can be formulated
mathematically as follows.\medskip
\end{subequations}
\begin{itemize}
\item Calder\'{o}n Problem\medskip
\end{itemize}

Let $\Omega$ be a bounded domain in $R^{n},n\geq2$, with a Lipschizian
boundary $\partial\Omega$, and $\gamma$ be a real bounded measurable function
in $\Omega$ with a positive lower bound. Consider the differential operator
\begin{equation}
L_{\gamma}(w)=\nabla\cdot(\gamma\nabla w), \label{eq6-1}%
\end{equation}
acting on $H^{1}(\Omega)$ and a quadratic form $Q_{\gamma}(\phi)$ for
functions in $H^{1}(R^{n})$%

\begin{equation}
Q_{\gamma}=\int_{\Omega}\gamma(\nabla w)^{2}dx,\ w\in H^{1}(R^{n}),\text{
}w|_{\partial\Omega}=\phi,
\end{equation}
where%

\[
L_{\gamma}(w)=0,\text{ \ in \ }\Omega.
\]

The problem is then to decide whether $\gamma$ is uniquely determined by the
operator $Q_{\gamma}$ and to calculate $\gamma$ in terms $Q_{\gamma}$ if
$\gamma$ is indeed determined by $Q_{\gamma}$. To put it in another way, we
need to verify the injectivity of the following map
\begin{equation}
\Phi:\gamma\rightarrow Q_{\gamma}. \label{eq6-4}%
\end{equation}

In the context of physical meaning, $Q_{\gamma}(\phi)$ represents the power
necessary to maintain an electrical potential $\phi$ on $\partial\Omega$.
Caldero\'{n} showed that the map $\Phi$ is analytic if $\gamma\in L^{\infty
}(\Omega)$ and $d\Phi|_{\gamma=const}$ is injective. He also proved that
$Q_{\gamma}$ determines $\gamma$ when $\gamma$ is sufficiently close to a constant.

Literally speaking, $Q_{\gamma}$ can be determined through boundary
measurements as we will see below. Then the problem is reduced to whether the
conductivity can be reconstructed through the surface electrode measurements.
Many authors have made contributions to the problem under various assumptions.
Kohn and Vogelius \cite{[Kohn1984]} provided a positive answer to the
determination of the conductivity which is $C^{\infty}$ in $\bar{\Omega}$ and
has all derivatives at the boundary. Sylvester and Uhlmann
\cite{[Sylvester1986]} proved uniqueness for $C^{2}$ conductivities in the
plane while Brown \cite{{[Russell1996]}} relaxed the regularity of the
conductivity to $3/2+\epsilon$ derivatives.

We may treat Calder\'{o}n Problem from another perspective. To be specific,
knowing $Q_{\gamma}(\phi)$ for each $\phi\in H^{1/2}(\Gamma)$ is equivalent to
the knowledge of \textquotedblleft Dirichlet-to-Neumann\textquotedblright%
\ data. In fact, by the Green's identity and assuming that the potential on
the boundary of the
\begin{equation}
\int_{\Omega\backslash\Omega_{0}}\left[  v\nabla\cdot(\gamma\nabla u)+\nabla
v\cdot\gamma\nabla u\right]  dx=\int_{\partial\left(  \Omega\backslash
\Omega_{0}\right)  }v\cdot\gamma\frac{\partial u}{\partial n}dS.\label{eq6-5}%
\end{equation}

If setting $v=u$ above, then we have
\begin{equation}
\int_{\Omega\backslash\Omega_{0}}\gamma|\nabla u|^{2}dx=\int_{\partial\left(
\Omega\backslash\Omega_{0}\right)  }\phi\cdot\gamma\frac{\partial u}{\partial
n}dS,\label{eq6-7}%
\end{equation}
where the left-hand side of (\ref{eq6-7}) is exactly $Q_{\gamma}(\phi)$ and
the right-hand side involves Neumann values if \ Dirichlet conditions are
given on the boundary. Thus, the \textbf{Calder\'{o}n Problem} can be restated
as whether $\gamma$ is uniquely determined by the \textquotedblleft
Dirichlet-to-Neumann\textquotedblright\ map on the boundary. The map tells us
how the boundary potential determines the current flux across the boundary
\cite{[Tararu2013]}. Also it is clear that a \textquotedblleft
Robin-to-Neumann\textquotedblright\ map (or voltage- to-current map) is
essentially equivalent to a \textquotedblleft
Dirichlet-to-Neumann\textquotedblright\ map through a close look at
(\ref{eq6-9-4}). Given a full Robin-to-Neumann map $R_{z_{l},\gamma}%
:\phi\rightarrow\gamma\nabla u|_{\partial\Omega}$, it uniquely determines
$z_{l}$ and thus equivalent to the Dirichlet-to-Neuman map. Under such
circumstances, uniqueness of solutions to the inverse conductivity problem was
proved by Astala and P\"{a}iv\"{a}rinta \cite{[Astala2006]} without any
regularity imposed on the boundary for a bounded measurable conductivity in
two dimensions. In three dimensions, Haberman and Tararu \cite{[Tararu2013]}
confirmed the answer for $C^{1}$ conductivities and Lipchitz conductivities
close to the identity.


The Robin-to-Neumann map is also called a voltage-to-current map. Under the
Complete Electrode model, to obtain the voltage-to-current map is equivalent
to solving the forward problem on the whole boundary given a known
conductivity. As discussed before, the traditional finite element method
or/and boundary element method may require dense mesh near the contacts of
skin and electrodes, resulting in large linear systems to be solved.
Therefore, more efficient schemes are desired without constraints from the
geometries of the domain and shapes of the electrodes.

One possible way to improve the efficiency would be to develop a probabilistic
estimator of the voltage-to-current map. The main advantage of the method lies
in the prevailing multicore computing. Maire and Simon \cite{[Maire2015]}
proposed a so-called partially reflecting random walk on spheres algorithm to
compute voltage-to-current map in a parallel manner, which is also efficient
when only solutions at only a few points are desired. For the Dirichlet
boundary problem, it is well known that killed Brownian motion is the
stochastic process that governs the differential operator. However, in the
mixed boundary situations, the partially reflecting Brownian motion comes into
play in preventing the path running out of the domain by either absorption or
instantaneous reflection. Simulation of absorption is much easier to take care
of comparing to that of reflection as the latter requires special techniques
like local finite difference discretization. Various schemes of first order or
second order schemes have been proposed and analyzed \cite{maire13}%
\cite{[Bossy2010]}. Maire and Simon also studied a similar approach involving
second order space discretization scheme. A variance reduction technique was
introduced as well to improve the efficiency and accuracy of the method.

In our work, we aim to find a probabilistic solution to the voltage-to-current
map by directly simulating the reflecting Brownian motion paths on the
boundary. More precisely, the calculation of the boundary local time is
treated explicitly in details and then integrated into the Feynman-Kac type
representation, and the voltages can be then obtained numerically on the boundary.


\section{Numerical scheme for forward problem and voltage-to-current map}

\begin{figure}[ptb]
\centering
\includegraphics[width=0.5\textwidth]{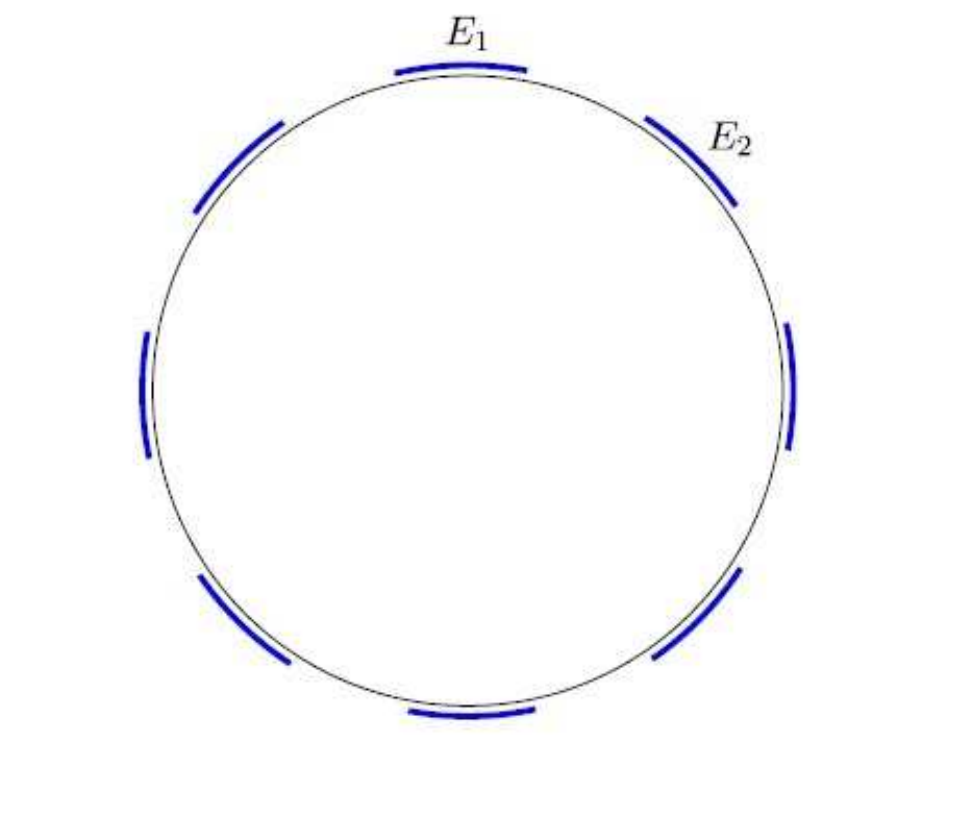}\caption{Limited
electrodes model on unit sphere}%
\label{fig:f6-1}%
\end{figure}

In medical applications, limited number of electrodes are attached to human
body to get surface measurements. We will illustrate our numerical method for
a model problem for a unit spherical object. In Fig. \ref{fig:f6-1}, eight
electrodes are superimposed on the boundary and the centers of the electrodes
all lie on the $y-z$ plane with radius 0.2.

Consider the conductivity equation (conductivity is taken to be 1 outside the
anomaly without loss of generality) with both Neumann and Robin boundary value
problem%
\begin{align}
\Delta u &  =0,\text{ \ \ in \ \ }\Omega_{0},\nonumber\\
z_{l}\nabla u\cdot n+u &  =\phi_{1}(x)=:\cos(4\theta)\text{ \ on }%
E_{l},l=1,\cdots,8,\label{mixedBC}\\
\frac{\partial u}{\partial n} &  =0,\text{ \ \ off }\cup_{l=1}^{8}%
E_{l},\nonumber\\
u &  =0\text{ \ \ \ \ \ on }\partial\Omega_{0}.
\end{align}
where $z_{l}$ is a constant between 0 and 1 and $n$ is the outward unit normal verctor.

In our previous work \cite{[Yijing2017]}, we described a method based on Monte
Carlo simulations to find potentials on boundaries through WOS sampling. The
same approach is employed to find voltages on the electrode patches given
mixed boundary conditions for the Laplace operator. From the Robin boundary
conditions on the electrode patches, the Neumann data are automatically known,
thus, the voltage-to-current map is obtained for the forward EIT\ problem.


\subsection{A path integral Monte Carlo Solution using Feynman-Kac formula}

First, let us review some preliminaries concerning Feynman-Kac formula for the
mixed boundary value problem in (\ref{mixedBC}), which lays the foundation of
our path integral Monte Carlo approach.


\subsubsection{Reflecting Brownian Motion and boundary local time}

Assume that $D$ is a domain with a $C^{1}$ boundary in $R^{3}$. A generalized
Skorohod problem is stated as follows: \newline

\begin{definition}
Let $f\in C([0,\infty),R^{3})$, a continuous function from $[0,\infty]$ to
$R^{3}$. A pair $(\xi_{t},L_{t})$ is a solution to the Skorohod equation
$S(f;D)$ if
\begin{enumerate}
\item $\xi$ is continuous in $\bar{D}$;
\item the local time $L(t)$ is a nondecreasing function which increases only
when $\xi\in\partial D$, namely,
\begin{equation}
L(t)=\int_{0}^{t}I_{\partial D}(\xi(s))L(ds); \label{eq6-33}%
\end{equation}
\item The Skorohod equation holds:
\begin{equation}
S(f;D):\qquad\ \xi(t)=f(t)-\frac{1}{2}\int_{0}^{t}n(\xi(s))L(ds),
\label{eq6-35}%
\end{equation}
where $n(x)$ denotes the outward unit normal vector at $x\in\partial D$.
\end{enumerate}
\end{definition}

The Skorohod problem was first studied in \cite{[R1]} by A.V. Skorohod in
addressing the construction of paths for diffusion processes inside domains
with boundaries, which experience the instantaneous reflection at the
boundaries. Skorohod presented the result in one dimension in the form of an
Ito integral and Hsu \cite{[Hsu1984]} later extended the concept to
$d$-dimensions ($d\geq2$).

In general, solvability of the Skorohod problem is closely related to the
smoothness of the domain $D$. For higher dimensions, the existence of
($\ref{eq6-35}$) is guranteed for $C^{1}$ domains while uniqueness can be
achieved for a $C^{2}$ domain by assuming the convexity for the domain
\cite{[15]}. Later, it was shown by Lions and Sznitman \cite{[16]} that the
constraints on $D$ can be relaxed to some locally convex properties.

Suppose that $f(t)$ is a standard Brownian motion (SBM) path starting at
$x\in\bar{D}$ and $(X_{t},L_{t})$ is the solution to the Skorohod problem
$S(f;D)$, then $X_{t}$ will be the standard reflecting Brownian motion (SRBM)
on $D$ starting at $x$. Because the transition probability density of the SRBM
satisfies the same parabolic differential equation as that for a BM, a sample
path of the SRBM can be simulated simply as that of \ the BM within the
domain. However, the zero Neumann boundary condition for the density of SRBM
implies that the path be pushed back at the boundary along the inward normal
direction whenever it attempts to cross the boundary.

The boundary local time $L_{t}$ is not an independent process but associated
with SRBM $X_{t}$ and defined by
\begin{equation}
L(t)\equiv\lim_{\epsilon\rightarrow0}\frac{\int_{0}^{t}I_{D_{\epsilon}}%
(X_{s})ds}{\epsilon}, \label{eq6-37}%
\end{equation}
where $D_{\epsilon}$ is a strip region of width $\epsilon$ containing
$\partial D$ and $D_{\epsilon}\subset\overline{D}$. Here $L_{t}$ is the local
time of $X_{t}$, a notion invented by P. L\'{e}vy \cite{[levy]}. This limit
exists both in $L^{2}$ and $P^{x}$-$a.s$. for any $x\in\overline{D}$.

It is obvious that $L_{t}$ measures the amount of time that the standard
reflecting Brownian motion $X_{t}$ spends in a vanishing neighborhood of the
boundary within the time period $[0,t]$. An interesting part of ($\ref{eq6-37}%
$) is that the set $\left\{  t\in R_{+}:X_{t}\in\partial D\right\}  $ has a
zero Lebesgue measure while the sojourn time of the set is nontrivial
\cite{[22]}. This concept is not just a mathematical one but also has physical
relevance in understanding the \textquotedblleft crossover exponent"
associated with \textquotedblleft renewal rate" in modern renewal theory
\cite{[17]}.


\subsubsection{Simulation of RBM and calculation of local time}

\noindent\textbf{Method of WOS for Brownian paths}

\medskip

Random walk on spheres (WOS) method was first proposed by M\"{u}ller
\cite{[muller]}, which can solve the Dirichlet problem for the Laplace
operator efficiently.

To illustrate the WOS method for the Dirichlet problem of Laplace equation,
with Dirichlet boundary conditions $\phi$. The solution can be rewritten in
terms of a measure $\mu_{D}^{x}$ defined on the boundary $\partial D$,
\begin{equation}
u(x)=E^{x}(\phi(X_{\tau_{D}}))=\int_{\partial D}\phi(y)d\mu_{D}^{x}%
,\label{eq6-38}%
\end{equation}
where $\mu_{D}^{x}$ is the harmonic measure\ defined by
\begin{equation}
\mu_{D}^{x}(F)=P^{x}\left\{  X_{\tau_{D}}\in F\right\}  ,F\subset\partial
D,x\in D.
\end{equation}
It can be shown easily that the harmonic measure is related to the Green's
function $g(y,x)$ for the domain with a homogeneous boundary condition
\cite{CKL}, i.e.,%
\begin{align*}
-\Delta g(x,y) &  =\delta(x-y),\text{ \ \ \ }x\in D,\\
g(x,y) &  =0,\text{ \ \ \ \ \ \ \ \ \ \ \ \ \ \ }x\in\partial D,
\end{align*}
as follows%
\begin{equation}
p(\mathbf{x},\mathbf{y})=-\frac{\partial g(x,y)}{\partial n_{y}}.
\end{equation}

\begin{figure}[ptb]
{\large \centering   \subfigure[WOS within the domain]{
\label{fig:subfig:a}     \includegraphics[width=0.37\textwidth]{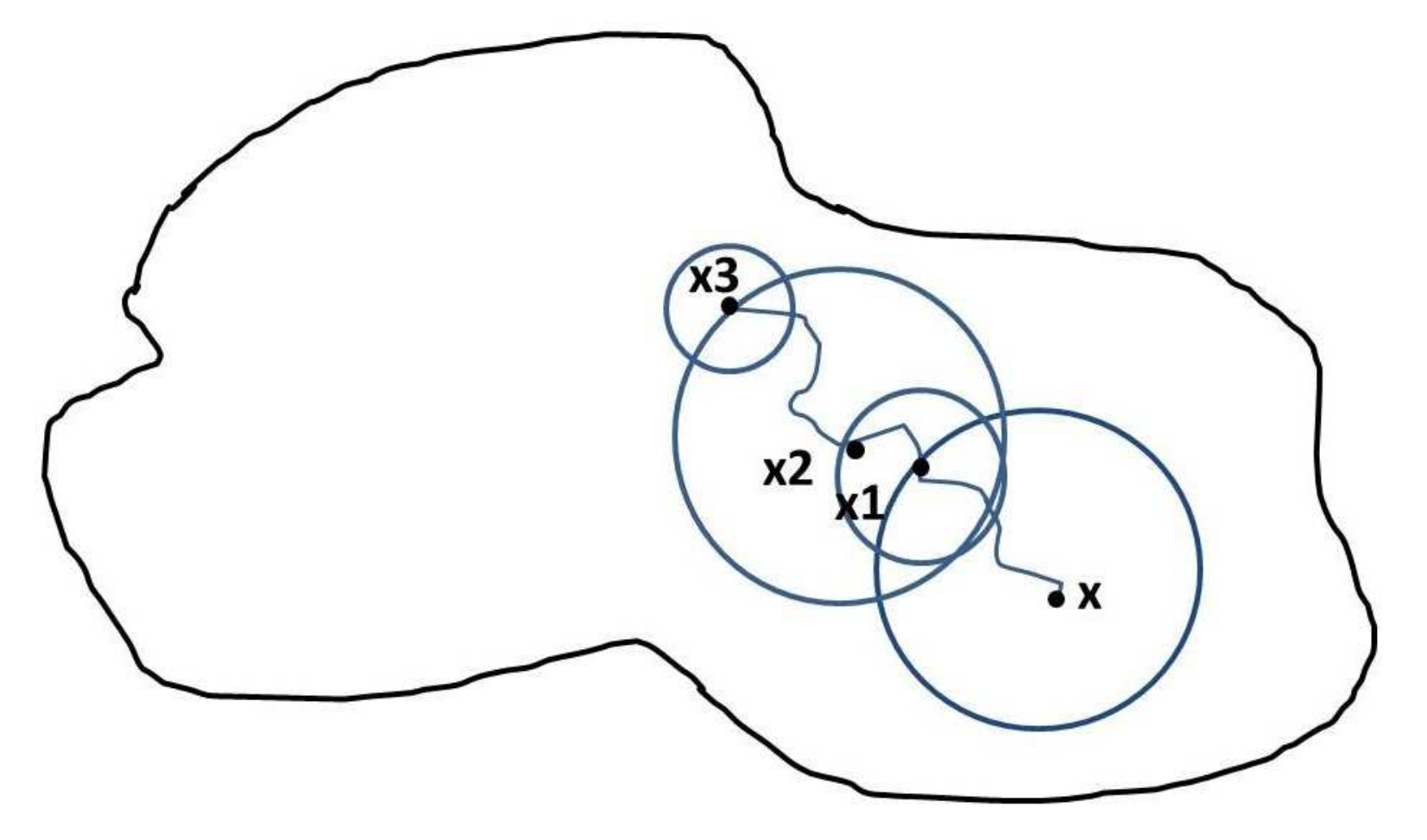}}
\hspace{1in}
\subfigure[WOS (with a maximal step size for each jump) within the domain]{
\label{fig:subfig:b}     \includegraphics[width=0.35\textwidth]{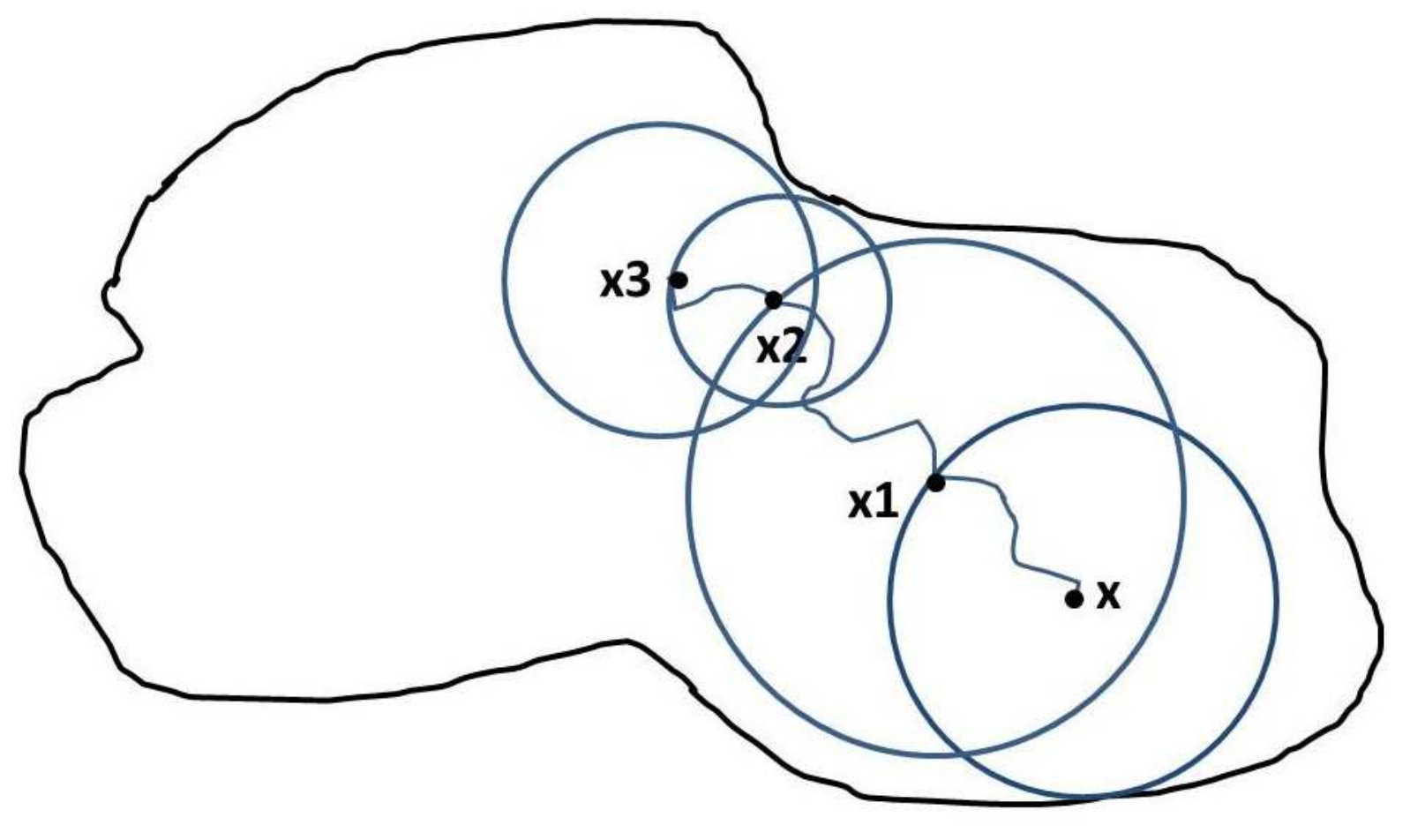}}
}\caption{Walk on Spheres method}%
\label{fig:subfig}%
\end{figure}

\begin{figure}[ptb]
\centering {\large \includegraphics[width=0.38\textwidth]{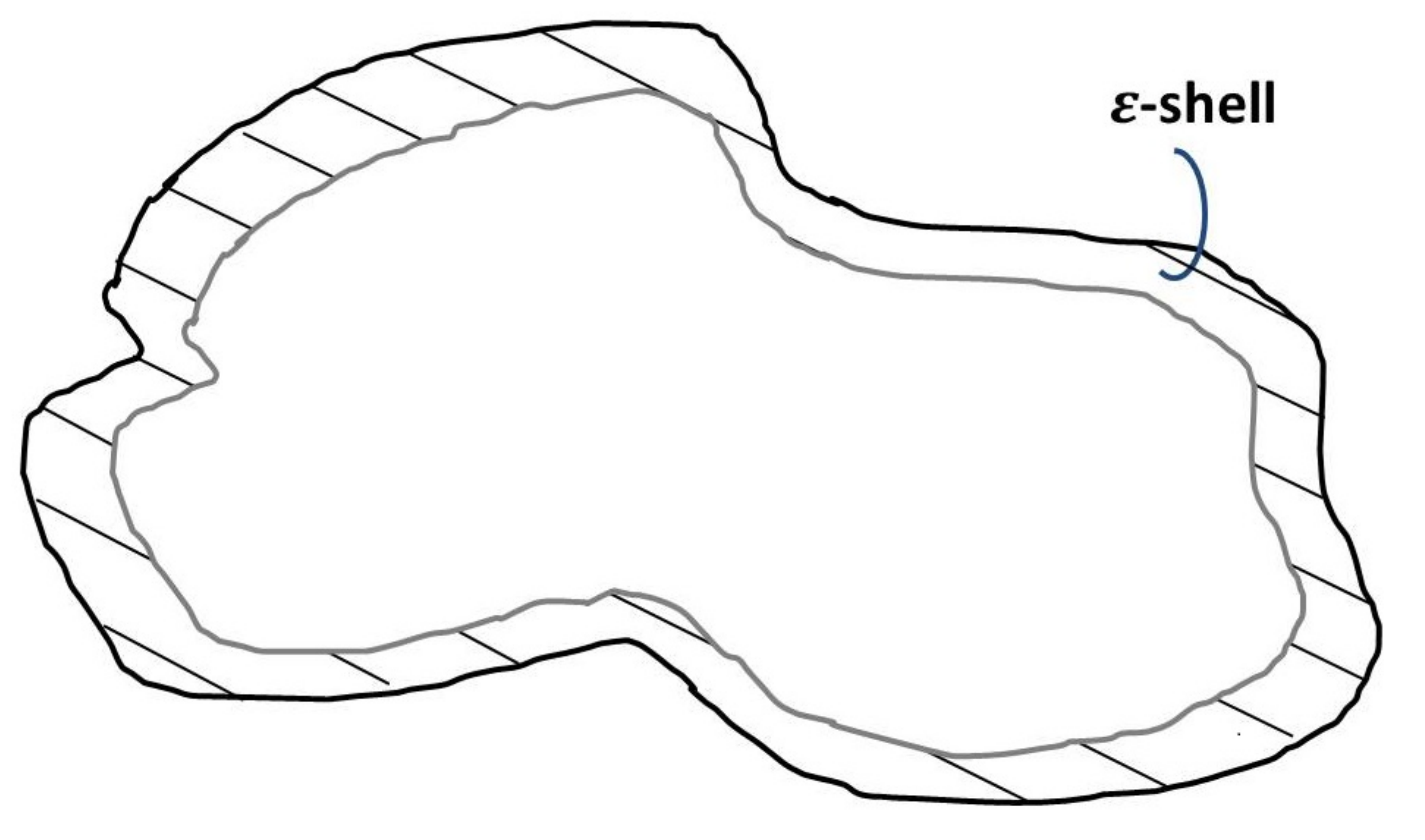}
}\caption{A $\epsilon$-region for a bounded domain in $R^{3}$}%
\label{fig:f6-11}%
\end{figure}

If the starting point $x$ of a Brownian motion is at the center of a ball, the
probability of the BM exiting a portion of the boundary of the ball will be
proportional to the portion's area. Therefore, sampling a Brownian path by
drawing balls within the domain can significantly reduce the path sampling
time. To be specific, given a starting point $x$ inside the domain $D$, we
simply draw a ball of largest possible radius fully contained in $D$ and then
the next location of the Brownian path on the surface of the ball can be
sampled, using a uniform distribution on the sphere, say at $x_{1}$. Treat
$x_{1}$ as the new starting point, draw a second ball fully contained in $D$,
make a jump from $x_{1}$ to $x_{2}$ on the surface of the second ball as
before. Repeat this procedure until the path hits a absorption $\epsilon
$-shell of the domain (see Fig. \ref{fig:f6-11}) \cite{[5]}. When this
happens, we assume that the path has hit the boundary $\partial D$ (see Fig.
\ref{fig:subfig:a} for an illustration).

Now we can define an estimator of (\ref{eq6-38}) by
\begin{equation}
u(x)\approx\frac{1}{N}\sum_{i=1}^{N}u(x_{i}),
\end{equation}
where $N$ is the number of Brownian paths sampled and $x_{i}$ is the first
hitting point of each path on the boundary. To speed up the WOS process,
maximum possible size of the sphere for each step would allow faster first
hitting on the boundary, see Fig. \ref{fig:subfig:b}.

\begin{figure}[ptb]
\centering {\large \includegraphics[width=0.5\textwidth]{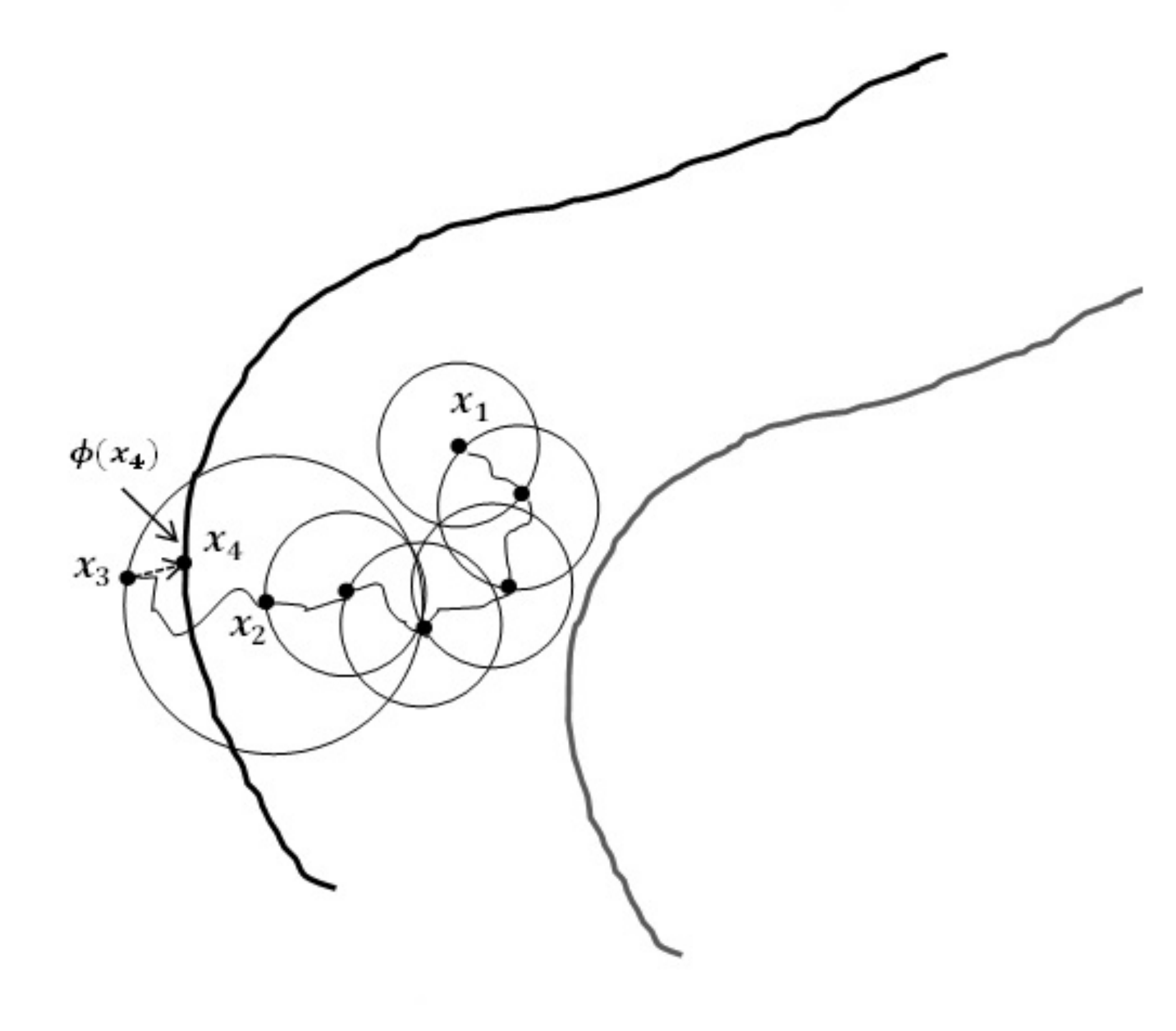}
}\caption{WOS in the $\epsilon$-region. At point $x_{1}$, BM path first hits
the $\epsilon$-region. By WOS with a prefixed radius $\Delta x$, the path
continues moving subsequently to $x_{2}$ where the distance to the boundary is
less than $\Delta x$. Enlarge the radius to $2\Delta x$, the path then have a
probability to run out of the domain to $x_{3}$. Pull back to the closest
point $x_{4}$ on the boundary, record $\phi(x_{4})$ and continue WOS-sampling
starting at $x_{4}$.}%
\label{fig:f6-13}%
\end{figure}

For the reflecting boundary, we will construct a strip region around the
boundary (see Fig. \ref{fig:f6-11}) and allow the process $X_{t}$ to move
according to the law of BM continuously. Before the path enters the strip
region, the radius of WOS is chosen to be of a maximum possible size within
the distance to the boundary. Once the particle is in the strip region, the
radius of the WOS\ sphere is then fixed at a constant $\Delta x$ (or $2\Delta
x$, see Fig. \ref{fig:f6-13}). With this approach, according to the definition
($\ref{eq6-37}$), the local time may be interpreted as
\begin{equation}
dL(t)\approx\frac{\int_{t_{j-1}}^{t_{j}}I_{D_{\epsilon}}(X_{s})ds}{\epsilon},
\label{eq6-39}%
\end{equation}
which is
\begin{equation}
dL(t)\approx\frac{\int_{t_{j-1}}^{t_{j}}I_{D_{\epsilon}}(X_{s})ds}{\epsilon
}=(n_{t_{j}}-n_{t_{j-1}})\frac{(\Delta x)^{2}}{3\epsilon}, \label{eq6-43}%
\end{equation}
given a prefixed constant $\Delta x$ in the strip region and $n_{t_{j}}$ be
the cumulative steps that path stays within the $\epsilon$-region from the
begining until time $t_{j}$ (see Remark below for definition). Notice that
only those steps where the path of $X_{t}$ remains in the $\epsilon$-region
will contribute to $n_{t_{j}}$ while the SRBM may lie out of the $\epsilon
$-region at other steps. More details can be found in \cite{[Yijing2017]}. One
may refer to Fig. \ref{fig:f6-13} for an illustration of the behavior of path
near the boundary.

\noindent\medskip

\textbf{Remark} Occupation time of SRBM $X_{t}$ in the numerator of
$(\ref{eq6-39})$ was calculated in terms of that of BM sampled by the walks on
spheres. Notice here that within the $\epsilon$-region, the radius of the WOS
may be $\Delta x$ or $2\Delta x$, which implies that the corresponding elapsed
time of one step for local time could be $(\Delta x)^{2}/3$ or $(2\Delta
x)^{2}/3$. The latter is four times bigger than the former. But if we absorb
the factor $4$ into $n_{t}$, $(\ref{eq6-43})$ still holds. In practical
implementation, we treat $n_{t}$ as a vector of entries of increasing value,
the increment of each component of $n_{t}$ over the previous one after each
step of WOS will be 0, 1 or 4, corresponding to the scenarios that $X_{t}$ is
out of the $\epsilon$-region, in the $\epsilon$-region while sampled on the
sphere of a radius $\Delta x$, or in the $\epsilon$-region while sampled on
the sphere of a radius $2\Delta x$, respectively.


\subsubsection{Feynamn-Kac formula}

We consider the mixed boundary value problem in the domain $\Omega
\backslash\Omega_{0}$ to the mixed problem%
\begin{align}
\Delta u  &  =0\text{ \ \ \ in }\Omega,\nonumber\\
\nabla u\cdot n+cu  &  =\phi_{1}(x)\text{ \ \ \ on \ }\Gamma_{1}=\cup
_{l=1}^{8}E_{l},\\
\frac{\partial u}{\partial n}  &  =\phi_{2}(x)\text{\ \ \ \ on \ }\Gamma
_{2}=\partial\Omega\backslash\Gamma_{1},\nonumber\\
u  &  =\phi_{3}(x)\text{ \ \ \ on \ }\Gamma_{3}=\partial\Omega_{0}.
\end{align}

The probabilistic solution for the boundary value problem above is given by
the well-known Feynman-Kac formula \cite{[8]}
\begin{align}
u_{Mix}(x)  &  =E^{x}\left\{  \int_{0}^{\infty}\hat{e_{c}}(t)\phi_{1}%
(X_{t})dL(t)\right\} \label{eq5-81}\\
&  +\frac{1}{2}E^{x}\left\{  \int_{0}^{\infty}\phi_{2}(X_{t})dL(t)\right\}
+E^{x}(\phi_{3}(X_{\tau_{\Gamma_{3}}})).\nonumber
\end{align}
where $X_{t}$ is the standard reflecting Brownian motion, $L(t)$ is the
corresponding local time and the Feynman-Kac functional $\hat{e_{c}%
}(t):=e^{\int_{0}^{t}c(X_{t})dL(t)}$, $\ \tau_{\Gamma_{3}}$ is the first time
a Brownian path originating from $\Omega\backslash\Omega_{0}$ hits the
boundary of $\partial\Omega_{0}=\Gamma_{3}$.

Feynman-Kac formula provides a local solution procedure to solve the partial
differential equations through stochastic processes. As a matter of fact, the
infinitesimal generator of the Laplace operator is Brownian motion which is
involved in the solution as well. For a general elliptical operator, It\'{o}
processes come into play.

The numerical approximation to $(\ref{eq5-81})$ will be
\begin{align}
\tilde{u}_{Mix}(x) &  =E^{x}\left\{  \int_{0}^{T}e^{\int_{0}^{t}c(X_{t}%
)dL(t)}\phi_{1}(X_{t})dL(t)\right\}  \label{eq5-83}\\
&  +\frac{1}{2}E^{x}\left\{  \int_{0}^{T}\phi_{2}(X_{t})dL(t)\right\}
+E^{x}(\phi_{3}(X_{\tau_{\Gamma_{3}}})),\nonumber
\end{align}
or
\begin{align}
\tilde{u}_{Mix}(x) &  =E^{x}\left\{  \sum_{j=0}^{NP}e^{\int_{0}^{t}c(X_{t_{j}%
})dL(t_{j})}\phi_{1}(X_{t_{j}})dL(t_{j})\right\}  \label{eq5-85}\\
&  +\frac{1}{2}E^{x}\left\{  \sum_{j=0}^{NP}\phi_{2}(X_{t_{j}})dL(t_{j}%
)\right\}  +E^{x}(\phi_{3}(X_{\tau_{\Gamma_{3}}})).\nonumber
\end{align}

Equivalently,
\begin{align}
\tilde{u}_{Mix}(x)  &  =E^{x}\left\{  \sum_{{j:X_{t_{j}}\in E}}^{NP}%
e^{\int_{0}^{t}c(X_{t_{j}})dL(t_{j})}\phi_{1}(X_{t_{j}})dL(t_{j})\right\}
\label{eq5-87}\\
&  +\frac{1}{2}E^{x}\left\{  \sum_{{j:X_{t_{j}}\in}\Gamma_{2}}^{NP}\phi
_{2}(X_{t_{j}})dL(t_{j})\right\}  +E^{x}(\phi_{3}(X_{\tau_{\Gamma_{3}}%
})).\nonumber
\end{align}

Recalling the approximation of the local time in ($\ref{eq6-43}$),
($\ref{eq5-87}$) can be modified as%
\begin{align}
\tilde{u}_{Mix}(x) &  =E^{x}\left\{  \sum_{{j}^{\prime}=0}^{NP}e^{\sum
_{k=0}^{j}c(X_{t_{k}})(n_{t_{k}}-n_{t_{k-1}})\frac{(\Delta x)^{2}}{3\epsilon}%
}\phi_{1}(X_{t_{j}})(n_{t_{j}}-n_{t_{j-1}})\frac{(\Delta x)^{2}}{3\epsilon
}\right\}  \nonumber\\
&  +\frac{1}{2}E^{x}\left\{  \sum_{{j}^{\prime}=0}^{NP}\phi_{2}(X_{t_{j}%
})(n_{t_{j}}-n_{t_{j-1}})\frac{(\Delta x)^{2}}{3\epsilon}\right\}
++E^{x}(\phi_{3}(X_{\tau_{\Gamma_{3}}})),\label{pimc}%
\end{align}
where $j^{\prime}$ denotes each step of the path and $j$ denotes the steps
where the path hits the boundary (Robin or Neumann).

In the context of complete electrode model, the second expectation is zero due
to the zero Neumann boundary. One may find more details in \cite{[Yijing2016]}.


\subsection{A deterministic solution with the boundary element method}

To provide reference solutions to the path integral MC method proposed above,
we will present a deterministic method based on boundary element method for
the mixed boundary value problem (\ref{mixedBC}).

\subsubsection{Graded boundary mesh}

Let $\Omega$ be the unit ball centered at the origin and $\partial
\Omega=\Gamma_{1}\cup\Gamma_{2}$, the collection of electrodes is $\Gamma
_{1}=\cup_{i=1}^{8}E_{i}$. Each electrode patch is assumed to have an equal
surface area.

\begin{figure}[ptb]
{\large \centering \subfigure[Construction of curved triangles in 2D]{
\includegraphics[width=0.5\textwidth]{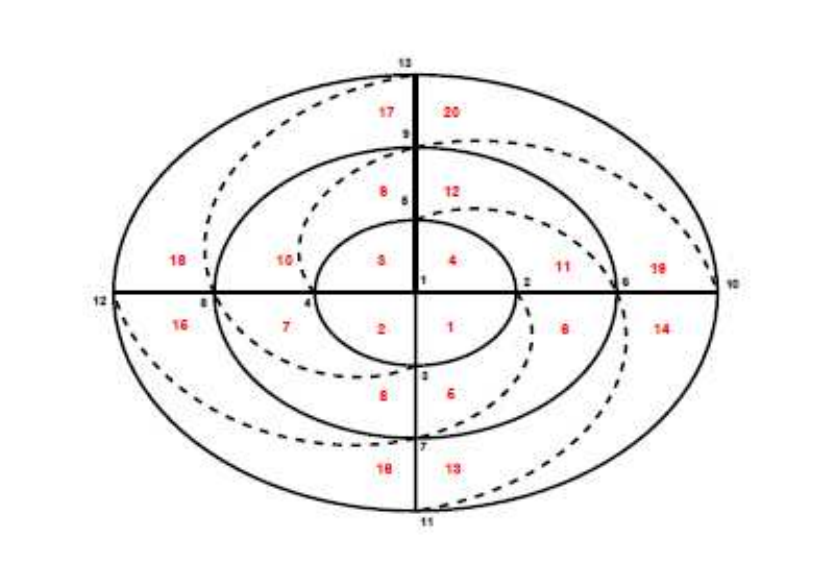}
\label{fig:f6-3a}} \subfigure[Mesh points of curved triangles]{
\includegraphics[width=0.5\textwidth]{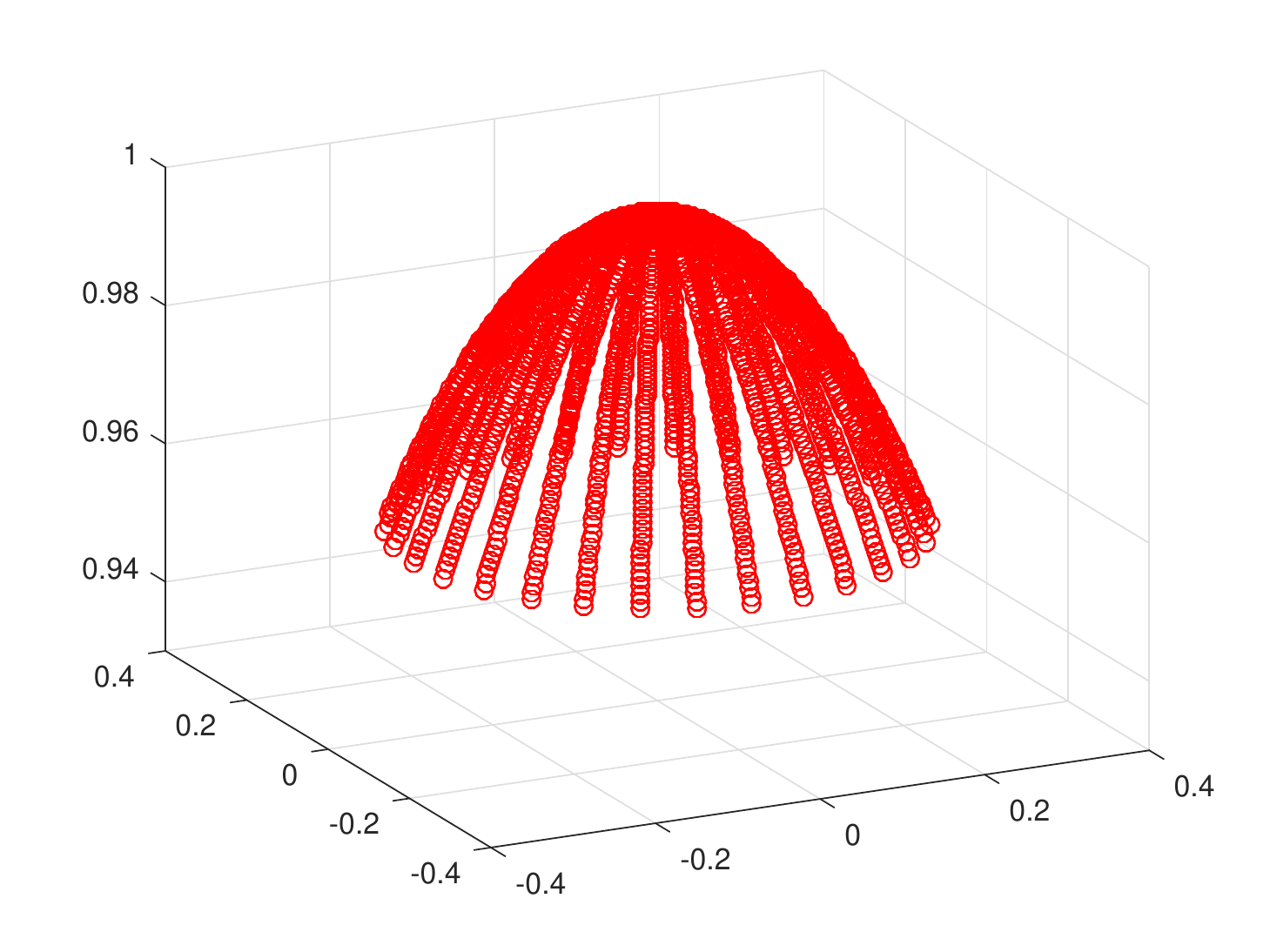}
\label{fig:f6-3b}
}}\caption{Curved triangles}%
\label{fig:f6-3}%
\end{figure}

\begin{figure}[ptb]
\centering
\includegraphics[width=0.5\textwidth]{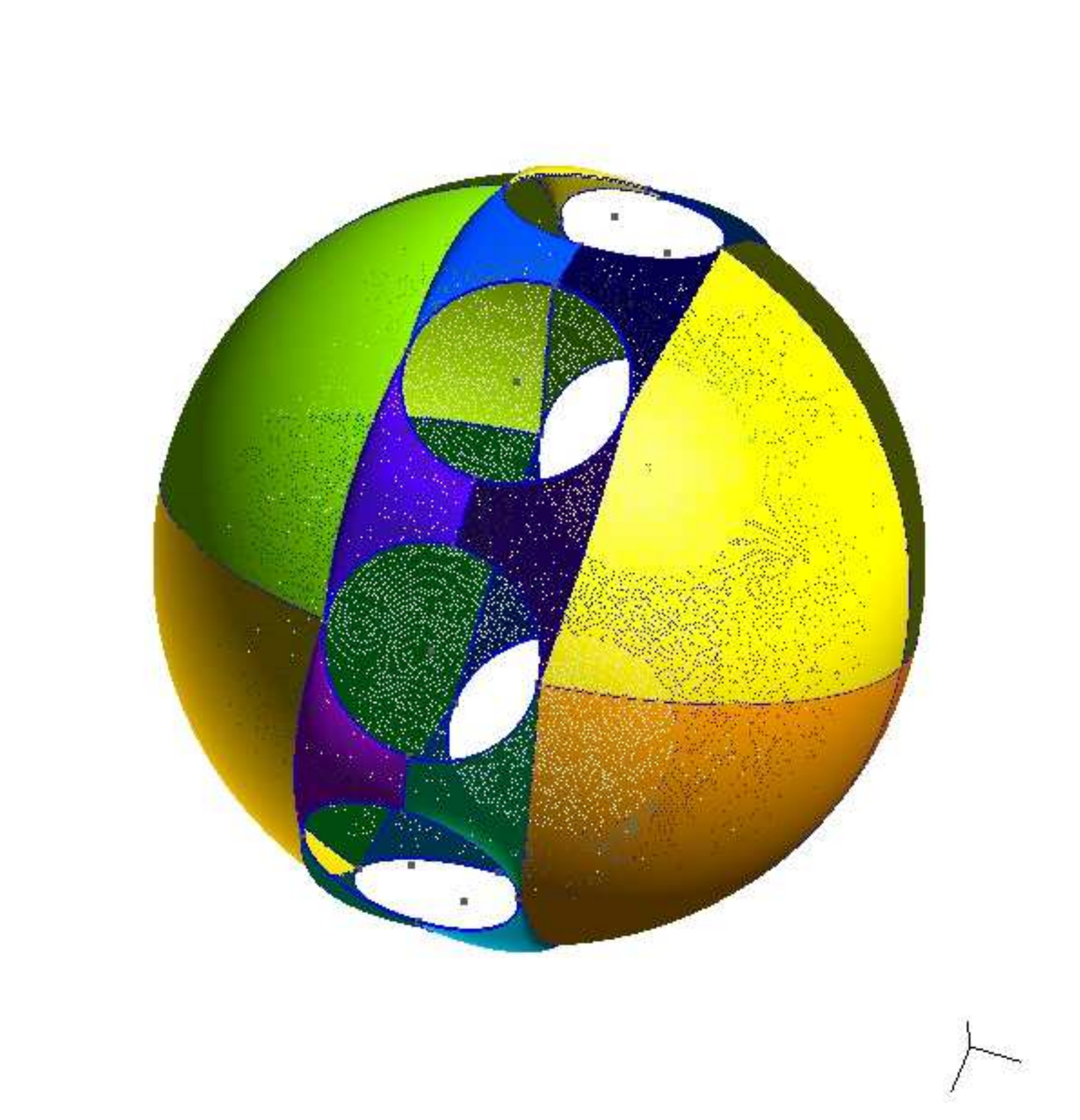}\caption{Boundary
mesh generated by GMSH on off-electrode patch.}%
\label{fig:f6-5}%
\end{figure}

Boundary mesh is constructed on the electrode patches and off-electrodes,
respectively as shown in Fig.\ref{fig:f6-3} and Fig.\ref{fig:f6-5}. For our
implementation, GMSH is used to generate an unstructured 2D mesh consisting of
flat triangles given a \textquotedblleft size field\textquotedblright\ while
on the electrode patches the meshes are structured in such a way that mesh
points are found at the intersection of division along the longitude and
altitude which gives a body-fitted mesh and the mapping from the elemental
triangle to the curved ones can be found in \cite{[cai2014]}. A global
boundary integral equation can then be set up based on the two meshes.

\begin{figure}[ptb]
\centering
\includegraphics[width=0.6\textwidth]{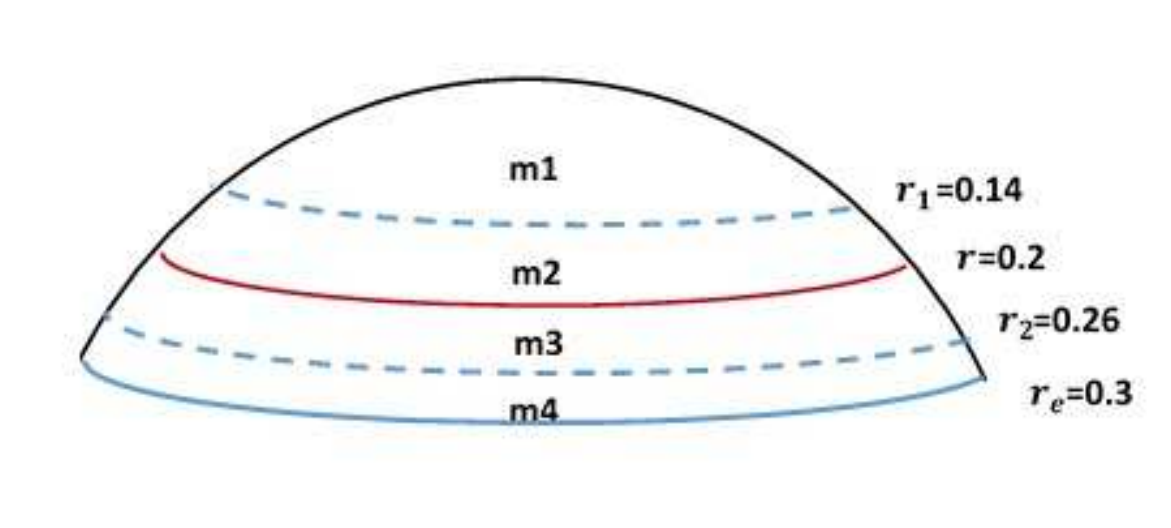}\caption{Graded mesh on
electrode patches.}%
\label{fig:f6-7}%
\end{figure}

\begin{figure}[ptb]
{\large \centering
\subfigure[View in 2D: Mesh points projected to x-y plane]{
\includegraphics[width=0.5\textwidth]{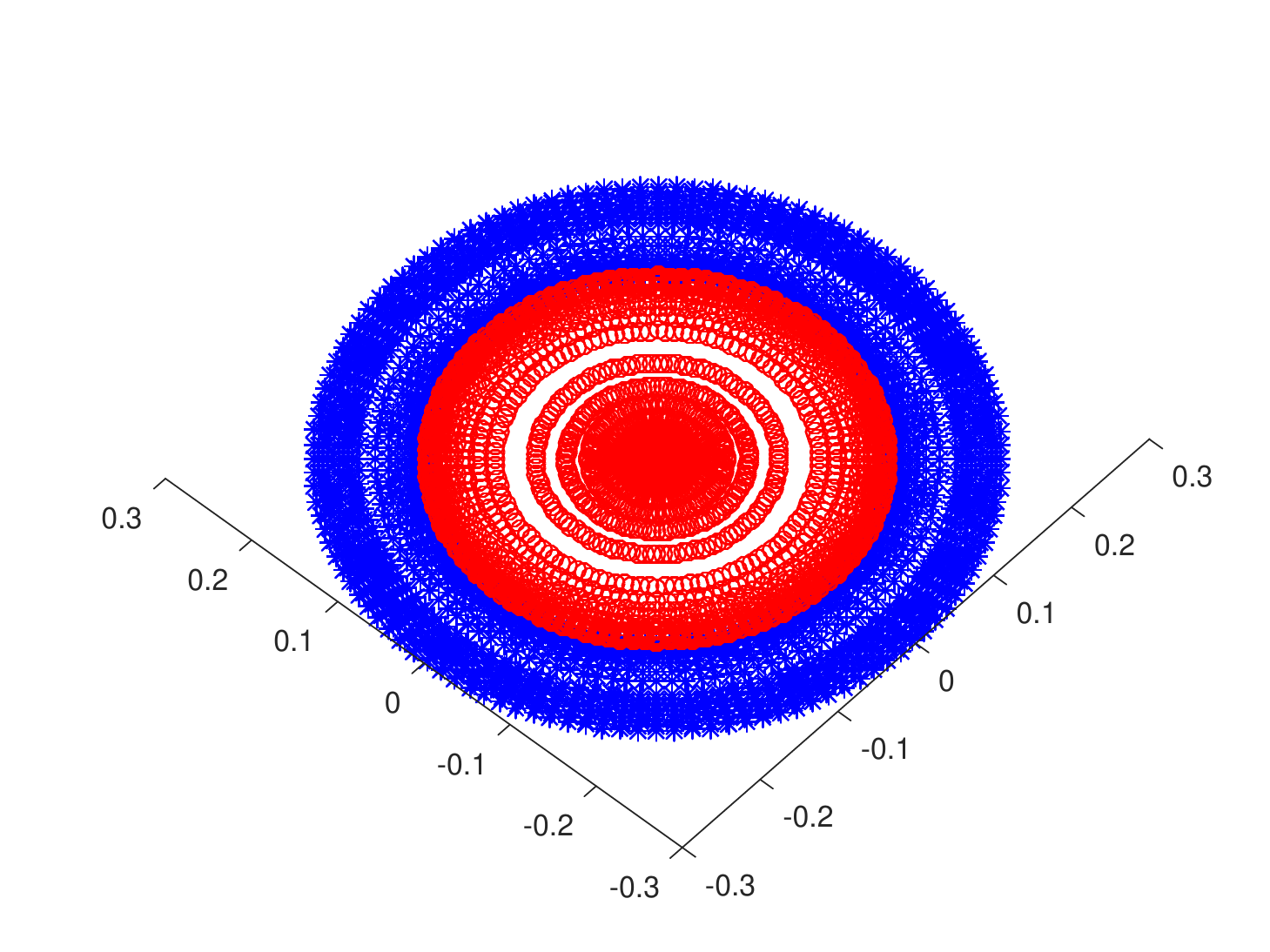}
\label{fig:f6-9a}} \subfigure[View in 3D]{
\includegraphics[width=0.5\textwidth]{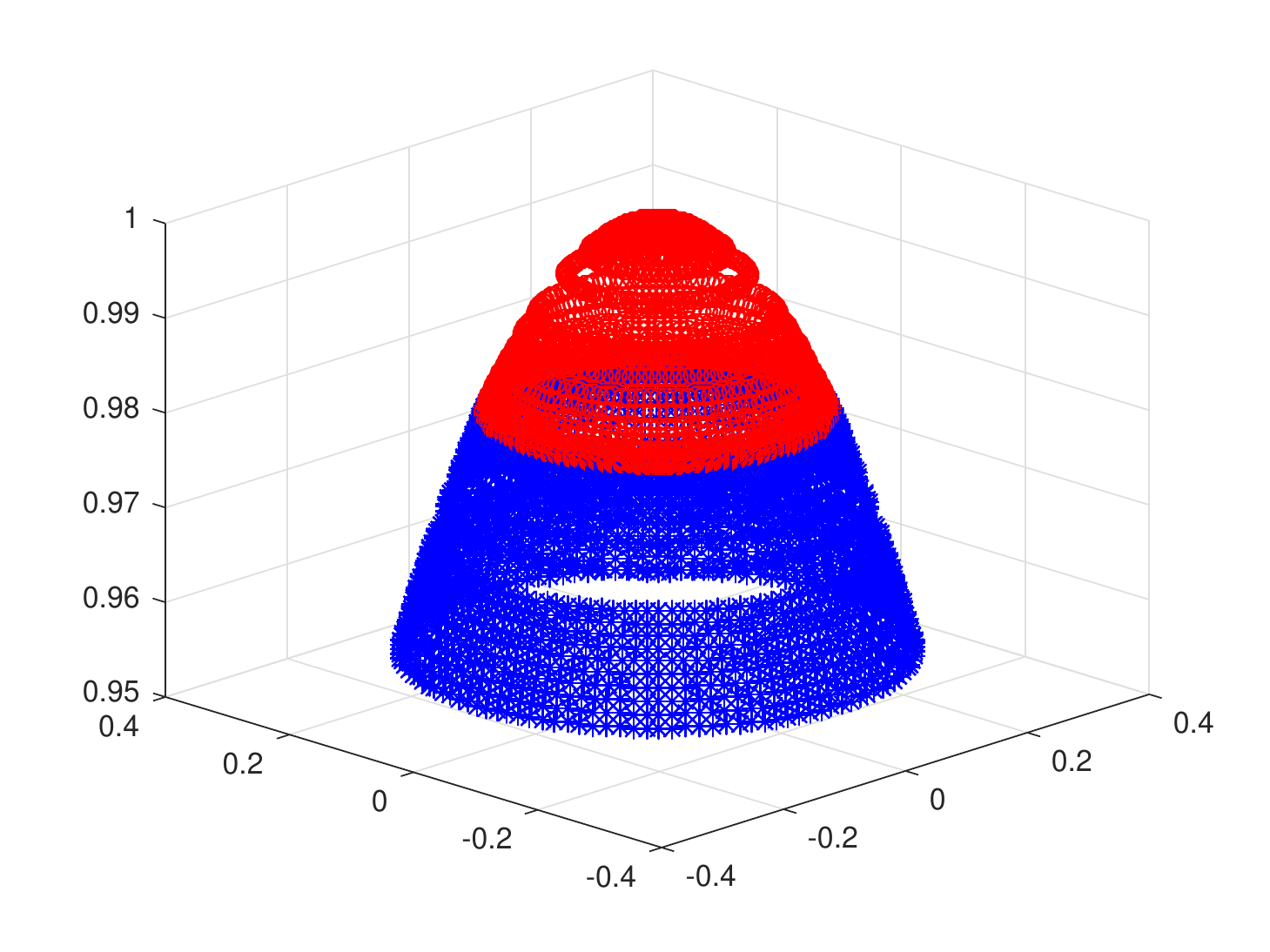}
\label{fig:f6-9b}
}}\caption{Mesh points on the north pole patch}%
\label{fig:f6-9}%
\end{figure}

Since the contact impedance $z_{l}$ varies from 0 to 1, we take $z_{l}=0.5$
for our numerical tests. A close look at the boundary conditions in
(\ref{mixedBC}) reveals that discontinuities at the rims of all the
electrodes. It is natural to enlarge the radius of mesh on $E_{i}$ so that we
may have an easy control over the mesh size for calculation. Assume the
enlarged radius to be $r_{e}=0.3$. Because of the discontinuity, we consider a
graded mesh on the enlarged surface by introducing a layered mesh structure,
as Fig. \ref{fig:f6-7} illustrates. There are four layers: the first ranges
from center to $r_{1}$, second from $r_{1}$ to $r$, third from $r$ to $r_{2}$
and fourth from $r_{2}$ to $r_{e}$. A dense mesh will be used around the rim
($r=0.2$) of the electrode, which implies that both 2nd layer and 3rd layer
should have a decreased mesh size towards $r=0.2$. Furthermore, a graded mesh
also discretizes the first layer while an evenly distributed mesh is used on
the fourth layer. And $m_{1},m_{2},m_{3}$ and $m_{4}$ are the number of
divisions along the altitude in each layer, respectively. Here we take
$m_{1}=20,m_{2}=16,m_{3}=16$ and $m_{4}=9$. The mesh size can be calculated
through $dx\cdot\alpha^{i},i=0,...,m_{j}-1(\alpha=3/4,j=1,2,3,4)$. The number
of divisions along longitude will be the same for each layer, i.e. $n=120$.
Fig. \ref{fig:f6-9} shows the realization of the graded mesh on the north pole
patch. The red points are the mesh points on the electrode (first two layers)
and blue ones are off-electrode points on the rest two layers. We can see
clearly the mesh on $E_{i},i=2,...,8$ can be constructed similarly or obtained
through rotation of that on $E_{1}$ along $x$-axis. Besides, the
\textquotedblleft size field\textquotedblright\ of GMSH is 0.012 on D which
yields 88383 mesh points and 175530 flat triangle elements as Fig.
\ref{fig:f6-5} shows.

\subsubsection{Boundary integral equation}

By using the second Green's identity, we have for any fixed $x\in
\Omega\backslash\Omega_{0}$,
\begin{align}
&  \int_{\Gamma}\left(  u(y)\Delta G(x,y)-G(x,y)\Delta u(y)\right)
dy\nonumber\\
&  =\int_{\Gamma}\frac{\partial G(x,y)}{\partial n_{y}}u(y)dS_{y}-\int
_{\Gamma}G(x,y)\frac{\partial u(y)}{\partial n_{y}}dS_{y},
\end{align}
where $\Gamma=\Gamma_{1}\cup\Gamma_{2}\cup\Gamma_{3}$ and $G(r,r^{\prime})$ is
a Green's function for the exterior domain of $\Omega_{0}$ with homogeuous
boundary condition on $\partial\Omega_{0}$, i.e.,%
\begin{align}
-\Delta G(x,y) &  =\delta(x-y)\ \ \ \text{in \ }x,\ \ y\in\text{\ }\Omega
_{0}^{c},\\
G(x,y)|_{x\in\partial\Omega_{0}} &  =0.\nonumber
\end{align}

As $x\rightarrow\Gamma_{1}\cup\Gamma_{2}$ from the interior, using the fact
that both $u$ and $G(x,y)$ vanish on the boundary $\Gamma_{3},$we obtain
\begin{subequations}
\label{eq6-27}%
\begin{align}
-\frac{1}{2}u(x) &  =\int_{\partial\Omega}\frac{\partial G(x,y)}{\partial
n_{y}}u(y)dS_{y}-\int_{\partial\Omega}G(x,y)\frac{\partial u(y)}{\partial
n_{y}}ds_{y},\label{eq6-27-1}\\
-\frac{1}{2}u(x) &  =\int_{\partial\Omega}\frac{\partial G(x,y)}{\partial
n_{y}}u(y)dS_{y}-\int_{\Gamma_{1}}G(x,y)\frac{\partial u(y)}{\partial n_{y}%
}ds_{y},\label{eq6-27-3}\\
-\frac{1}{2}u(x) &  =\int_{\partial\Omega}\frac{\partial G(x,y)}{\partial
n_{y}}u(y)dS_{y}-\int_{\Gamma_{1}}G(x,y)\left(  \phi_{1}(y)-u(y)\right)
/z_{l}ds_{y},\label{eq6-27-5}\\
-\frac{1}{2}u(x) &  =\int_{\partial\Omega}\frac{\partial G(x,y)}{\partial
n_{y}}u(y)dS_{y}-\frac{1}{z_{l}}\int_{\Gamma_{1}}G(x,y)\phi_{1}(y)dS_{y}%
\label{eq6-27-7}\\
&  +\frac{1}{z_{l}}\int_{\Gamma_{1}}G(x,y)u(y)dS_{y},\nonumber
\end{align}

or
\end{subequations}
\begin{align}
&  \frac{1}{2}u(x)+\int_{\partial\Omega}\frac{\partial G(x,y)}{\partial n_{y}%
}u(y)dS_{y}+\label{eq6-29}\\
&  \frac{c}{z_{l}}\int_{\Gamma_{1}}G(x,y)u(y)dS_{y}=\frac{1}{z_{l}}%
\int_{\Gamma_{1}}G(x,y)\phi_{1}(y)ds_{y},\nonumber
\end{align}
where the Robin condition in (\ref{mixedBC}) is used to deduce from
(\ref{eq6-27-1}) to (\ref{eq6-27-3}). Therefore we obtain a global boundary
integral equation for $u$. With the different meshes on $D$ and $E$, we let
$x$ sweep over all the mesh points where (\ref{eq6-29}) is imposed, we have a
linear system for solving $u(x)$ on $\partial\Omega$. As a result, reference
potentials on the boundary can be obtained to validate the Monte Carlo
simulations. Meanwhile, any reference potentials inside the domain can be
achieved as well through
\begin{equation}
u(x)=\int_{\partial\Omega}G(x,y)\frac{\partial u(y)}{\partial n_{y}}%
dS_{y}-\int_{\partial\Omega}\frac{\partial G(x,y)}{\partial n_{y}}%
u(y)dS_{y}.\label{eq6-31}%
\end{equation}
where $\frac{\partial u(y)}{\partial n_{y}}$, the Neumann values on the
boundary, are automatically known once the reference potentials are found on
the electrode patches from the Robin conditions.


\subsubsection{Reference current on electrodes}

With the preparations in the previous sections, the results are shown in terms
of current on each electrode according to
\begin{equation}
J_{l}^{ref}=\frac{1}{|E_{l}|}\int_{E_{l}}\frac{\partial u}{\partial
n}\bigg\rvert_{\partial D}d\sigma(x), \label{eq6-32}%
\end{equation}
\begin{table}[ptb]
\caption{Reference currents on each electrode. \ }%
\label{tab:table6-1}
\begin{center}%
\begin{tabular}
[c]{|c|c|c|c|c|}\hline
Current & $E_{1}$ & $E_{2}$ & $E_{3}$ & $E_{4}$\\\hline
$J^{ref}$ & 1.3377346024 & -1.3960453685 & 1.4543557058 &
-1.3960453502\\\hline
Current & $E_{5}$ & $E_{6}$ & $E_{7}$ & $E_{8}$\\\hline
$J^{ref}$ & 1.3377346471 & -1.3960453364 & 1.4543557459 &
-1.3960453565\\\hline
\end{tabular}
\end{center}
\end{table}

A direct summation of $J_{l}^{ref}$ over eight electrodes yields the whole
current to be \textbf{ -7.10e-7} which is close to \textbf{0} as the
conservation of charges condition suggests. Meanwhile the electrode currents
show symmetric patterns with respect to both $y$ and $z$ axis, consistent with
the system design.


\subsection{Validation of path integral MC numerical results}

All Monte Carlo simulation is based on the algorithm in (\ref{pimc}). The
number of Monte Carlo simulations $N$ is $2\times10^{5}$ for all the mesh
points on the electrodes while the length of path $NP$ ranges from 1900 to
3000 on different patches. From Table \ref{tab:table6-3}, the absolute errors
between the numerical approximations and reference currents remain at a low
level not higher than 1\% .

Walk-on-sphere method is used to simulate the reflecting Brownian motion
within the spherical object. Feynman-Kac formula shows an infinite length of
path while from the numerical perspective it is truncated to be a finite
number $NP$ in such a way that $NP$ is increased until the desired accuracy is
achieved and it is lowered when $NP$ is increased further through our test.

With the numerical approximations of potentials on the boundary, the Neumann
values are automatically known and thus a full Robin-to-Neumann map has been
achieved over the whole boundary by conducting the same procedure at all the
mesh points on different electrode patches. Therefore, we've found an
effective way to compute the voltage-to-current map without resorting to the
conventional deterministic methods which involves constructing 3-D mesh in the
domain and solution of a large global matrix. It has also advantages when the
geometry of the domain in 3-D is complicated.

\begin{table}[ptb]
\caption{Numerical approximations of the current on eight patches}%
\label{tab:table6-3}
\begin{center}%
\begin{tabular}
[c]{|c|c|c||c|c|c|}\hline
$J_{1}^{ref}$ & $\tilde{J_{1}}$ & $|\text{Err}|$ & $J_{2}^{ref}$ &
$\tilde{J_{2}}$ & $|\text{Err}|$\\\hline
1.33773 & 1.33756 & 0.01\% & -1.39604 & -1.38565 & \ 0.74\%\\\hline
$J_{3}^{ref}$ & $\tilde{J_{3}}$ & $|\text{Err}|$ & $J_{4}^{ref}$ &
$\tilde{J_{4}}$ & $|\text{Err}|$\\\hline
1.45435 & 1.45521 & \ 0.06\% & -1.39604 & -1.39432 & \ 0.12\%\\\hline
$J_{5}^{ref}$ & $\tilde{J_{5}}$ & $|\text{Err}|$ & $J_{6}^{ref}$ &
$\tilde{J_{6}}$ & $|\text{Err}|$\\\hline
1.33773 & 1.33594 & \ 0.13\% & -1.39604 & -1.39717 & \ 0.08\%\\\hline
$J_{7}^{ref}$ & $\tilde{J_{7}}$ & $|\text{Err}|$ & $J_{8}^{ref}$ &
$\tilde{J_{8}}$ & $|\text{Err}|$\\\hline
1.45435 & 1.45770 & 0.23\% & -1.39604 & -1.38508 & 0.79\%\\\hline
\end{tabular}
\end{center}
\end{table}


\section{Conclusions and future work}

This paper presents a path integral Monte Carlo method to solve the forward
problem of EIT and the voltage-to-current map is acquired as needed for the
iterative algorithm of the inverse EIT problem where the forward problem needs
to be solved first at each iteration. The method takes advantage of the
parallel computing capability in modern multicore computers and solutions can
be calculated simultaneously at different electrodes independently.

In our calculations, the contact impedance $z_{l}$ is taken to be 0.5 while it
can be very close to zero. In that case, the voltages will jump drastically
and then a dense mesh must be placed at the contacts otherwise the calculation
of the reference voltages may be undermined. However, the PIMC will not be
affected as each point is independently calculated.

Meanwhile, local boundary integral equations may be considered as an
alternative approach to the voltage-to-current map. To be more specific, when
the reference potentials are obtained through the global boundary integral
equation, then potentials at any interior points are known due to
(\ref{eq6-31}). Then, we may construct local boundary integral equations of
Neumann values on each electrode.

The illness of the inverse problem requires high accuracy of forward modelling
otherwise it may lead to very large fluctuations in reconstructions of
conductivity. Our method is accurate to some extent but needs more development
in convergence and variance reduction.

\section*{Acknowledgement}

W.C. acknowledges the support of the National Science
Foundation (DMS-1764187) for the work in this paper.

\end{document}